\documentclass[11pt]{article}
\usepackage{epsfig}
\usepackage{amsmath,amsfonts,amsthm}
\newcommand{\integers}{\mathbb{Z}}
\newcommand{\figswitch}[2]{\epsfig{#1,#2} }                   
\theoremstyle{plain}
\newtheorem{thm}{Theorem}
\newtheorem{lem}{Lemma}
\newtheorem*{cor}{Corollary}
\newtheorem{prop}{Proposition}
\theoremstyle{definition}
\newtheorem{defin}{Definition}

\newcommand{\evec}[3]{ 
\stackrel{#1}{\varphi}{\!\!}^{\mbox{\tiny #2}}_{{#3} } } 
\newcommand{\Evec}[3] { 
\stackrel{#1}{\Phi}{\!\!}^{\mbox{\tiny #2}}_{{#3} } }
\newcommand{\evecb}[3]{ 
\stackrel{#1}{\varphi}{\!\!}^{\mbox{\tiny #2}*}_{{#3} } } 
\newcommand{\Evecb}[3] { 
\stackrel{#1}{\Phi}{\!\!}^{\mbox{\tiny #2}*}_{{#3} } }

\newcommand{\ZNWW}{{\stackrel{=}{Z}}^{\rm{\cal{N}}}}

\newcommand{\ZSW}{\bar{Z}^{\rm{\cal{S}}}}
\newcommand{\ZSWW}{{\stackrel{=}{Z}}^{\rm{\cal{S}}}}
\newcommand{\ZNW}{\bar{Z}^{\rm{\cal{N}}}}

\newcommand{\cU}[2]{\,\stackrel{{\scriptscriptstyle \, #1}}{\cal U}_{#2}}  
\newcommand{\cS}[2]{\,\stackrel{{\scriptscriptstyle \, #1}}{\cal S}_{#2}}
\newcommand{\Teo}{\stackrel{eo}{\mathbf{T}}_{1}  }
\newcommand{\Toe}{\stackrel{oe}{\mathbf{T}}_{1}}
\newcommand{\TNeo}{{\stackrel{eo}{\mathbf{T}}_{N}}  }
\newcommand{\TNoe}{{\stackrel{oe}{\mathbf{T}}_{N} } }

\newcommand{\Te}{\stackrel{e\cdot e}{\mathbf{T}}_{1}} 
\newcommand{\To}{\stackrel{o\cdot o}{\mathbf{T}}_{1}} 
\newcommand{\TNe}{\stackrel{e\cdot e}{\mathbf{T}}_{N}} 
\newcommand{\TNo}{\stackrel{o\cdot o}{\mathbf{T}}_{N}}

\newcommand{\cbK}{{\cal K}_{N}}
\newcommand{\cK}{{\cal K}}

\newcommand{\Lam}{\Lambda} 
\def\a{\alpha} 
\def\s{\sigma}

\newcommand{\e}{\epsilon}

\def\d{\delta} 
\def\e{\epsilon}

\def\l{\lambda}

\newcommand{\xf}{y^{f}}

\newcommand{\bx}{{\bf y}}
\newcommand{\bxi}{{\bf y}^{i}}
\newcommand{\bxf}{{\bf y}^{f}}
\def\bk{{\bf k}}
\renewcommand{\xi}{{y}^{i}}

\newcommand{\by}{{\bf y}}
\newcommand{\byi}{{\bf y}^{i}}
\newcommand{\byf}{{\bf y}^{f}}

\def\bx{{\bf y}}
\def\bk{{\bf k}}

\newcommand{\pto}{\!\rightarrow\!}

\newcommand{\suchthat}{\,|\,}

\def\half{\frac{1}{2}}

\begin{document}
\pagenumbering{roman} \setcounter{page}{0} 
\title{From the Bethe Ansatz to the Gessel-Viennot Theorem}
\author{R. Brak\dag, J. W. Essam\ddag \, and A. L. Owczarek\dag
        \thanks{{\tt {\rm email:} brak@maths.mu.oz.au,
j.essam@vms.rhbnc.ac.uk, aleks@ms.unimelb.edu.au}} \\
         \dag Department of Mathematics and Statistics,\\
         The University of Melbourne,\\
         Parkville, Victoria 3052, Australia\\
         \ddag Department of Mathematics,\\
         Royal Holloway, University of London,\\
         Egham, Surrey TW20 0EX, England.}
\date{
\begin{center}
	23 February, 1999 
\end{center}
}

\maketitle 
 
\begin{abstract} 
We state and prove several theorems that demonstrate how the
coordinate Bethe Ansatz for the eigenvectors of suitable transfer
matrices of a generalised inhomogeneous five-vertex model on the
square lattice, given certain conditions hold, is equivalent to the
Gessel-Viennot determinant for the number of configurations of $N$
non-intersecting directed lattice paths, or vicious walkers, with
various boundary conditions. Our theorems are sufficiently general to
allow generalisation to any regular planar lattice.
\vspace{1cm} 
 

\noindent{\bf Key words:} Vicious walkers, Lattice Paths,
Gessel-Viennot Theorem, Bethe Ansatz, Transfer Matrix method. 

\end{abstract} 
\vfill

\newpage
\pagenumbering{arabic}

\section{Introduction} 
\setcounter{page}{1}
The problem of non-intersecting paths or vicious walkers has been
studied by the statistical mechanics community who have been
interested in them as simple models of various polymer and other
physical systems and independently by the combinatorics community who
have been interested in them in connection with binomial determinants.
In statistical mechanics they are generally known as vicious walkers,
a term coined by Fisher \cite{fisher84}, who studied the continuous
version of the model.  Various cases of the lattice problem were later
studied by Forrester
\cite{forrester89,forrester89b,forrester90}.  Independently in the area of
combinatorics the problem of non-intersecting paths was solved by a
very general theorem of Gessel and Viennot \cite{gessel85,gessel89},
following the work of Lindstr\"{o}m \cite{lindstrom73}, and Karlin and
McGregor \cite{karlin59b,karlin59}.  All these studies express the
number of configurations as the value of a
determinant. Non-intersecting walks arose in yet another context, that
of vertex models in statistical mechanics, where it was noticed that
if the vertices of the six-vertex model are drawn in a particular way
they could be interpreted as lattice paths
\cite{wu68a,guttmann98a}.  If one of the vertices had weight zero,
giving a five-vertex model, the resulting paths were
non-intersecting. The vertex models are traditionally solved by
expressing the partition function (a generating function) in terms of
transfer matrices.  The partition function is then evaluated by either
of two very powerful techniques, that of commuting transfer matrices
\cite{baxter82} or by direct diagonalisation of the transfer matrices
using the coordinate Bethe Ansatz \cite{bethe31,lieb67b}. In this
paper we bring the independent results of the two communities together
for the case of $N$ non-intersecting paths.  We will show that the
Bethe Ansatz (from statistical mechanics) and the Gessel-Viennot
Theorem (from combinatorics) are essentially equivalent for a fairly
general problem on the square lattice. The theorems proved should be
easily generalisable to other planar lattices. The connection between the
six-vertex model and non-intersecting path problems, from which this
correspondence stems, has been recently discussed in general terms
\cite{guttmann98a}. Here we shall consider a model equivalent to a
generalised five-vertex model on the square lattice where the
(Boltzmann) weights associated with walk edges are inhomogeneous in
one direction.

\section{The model}

A lattice path or walk in this paper is a walk on a square lattice
rotated $45^{\circ}$ which has steps in only the north-east or
south-east directions, and with sites labelled $(m,y)$ (see figure 1). A set
of walks is \emph{non-intersecting} if they have no sites in common.
We are concerned with enumerating the number of configurations of $N$
non-intersecting walks, starting and ending at given positions, in
various geometries: 1) walks in a plane without boundaries; 2) walks
which are confined to the upper half plane; and 3) walks which are
confined to a strip of a given width, $L$.  More generally, one may be
interested in interacting cases where the walks nearest the boundaries
are attracted or repulsed by contact interactions: combinatorially
this requires knowledge of the number of walks with particular numbers
of contacts with each of the boundaries.  In this paper we shall focus
on case 3 since the other 2 cases can be easily obtained from this
case as limits.

To more easily describe our model we require the following sub-domains
of $\integers^{N}$
\begin{subequations}
	\label{eq:doms}
\begin{align}
	\cS{o}{L}
&=\{y \suchthat 1\le y\le L,
	y\in\integers \text{ and $y$ odd}\},
	\label{eq:stripnoo}\\
	\cS{e}{L}
 &=\{y \suchthat 0\le y\le L,
	y\in\integers \text{ and $y$ even}\},
	\label{eq:stripnoe}\\
	\cS{}{L}
&=\{y \suchthat 0\le y\le L,y\in\integers\},
	\label{eq:stripno}\\
	\cU{o}{L}
&=  \{(y_{1},\ldots,y_{N})\,| \, 1\le y_{1} 
		< \ldots <y_{N}\le L, \text{ $y_{i}\in\cS{o}{L}$}\} \label{eq:doms3}\\
	\cU{e}{L} 
&= \{(y_{1},\ldots,y_{N})\,| \, 0\le y_{1} < 
		\ldots <y_{N}\le L, \text{ $y_{i}\in\cS{e}{L}$}\}
		\label{eq:doms4}\\
	\cU{}{L} 
&= \{(y_{1},\ldots,y_{N})\,| \, 0\le y_{1} < 
		\ldots <y_{N}\le L, \text{ $y_{i}\in\cS{}{L}$}\} 
		\label{eq:doms5}
\end{align}
\end{subequations}
We will use $\cU{p}{L}$ to denote $\cU{o}{L}$ or $\cU{e}{L}$.
Let $N$ non-intersecting walks, confined to a strip of width $L$,
start at $y$-coordinates $\bxi=(\xi_{1},\ldots,\xi_{N})\in\,
\cU{p}{L}$ in column $m=0$ of the lattice sites
and terminate after $t$ steps at $y$-coordinates
$\bxf=(\xf_{1},\ldots,\xf_{N})\in\, \cU{p'}{L}$ in the $t^{th}$
column. If $t$ is even then $p'=p$ else $p'=\bar p$, where $\bar{p}$
is the opposite parity to $p$. We will only consider the case that $L$
is odd so that $|\!\cU{e}{L}\!|=|\!\cU{o}{L}\!|={{\half(L+1)}\choose N}$.  (If
$L$ is even a null space enters the subsequent analysis of the
transfer matrices leading to a distracting complication.) We are
considering paths such that a) if $(m-1,y)$ is the position of a path
in column $m-1$ the only possible positions for that path in column
$m$ are $(m,y^{\prime})$ with $y^{\prime}=y\pm 1$ and $0\le
y^{\prime}\le L$ and b) the non-intersection is defined through the
constraint that if there are $N$ sites occupied at $m=0$ then in each
column of sites ($0\leq m\leq t$) there are exactly $N$ occupied
sites. We generalise the walk problem associated with the five-vertex
problem \cite{wu68a,guttmann98a} by assigning a weight $w(y,y^{\prime})$ to
the lattice edge from site $(m-1,y)$ to $(m,y^{\prime})$ with
$y^{\prime}=y\pm 1$ (see figure 1). Notice that, since
$w(y,y^{\prime})$ is assumed independent of the column index $m$, due
to the square lattice structure the weights are periodic in the $m$
direction with period two: Note if $y\in\cS{p}{L}$ then
$y'\in\cS{\bar{p}}{L}$, and in general $w(y,y')\ne w(y',y)$. For the sake
of generality we also associate an arbitrary weight $v(y^i)$ with each
of the sites occupied at $m=0$. The weight associated with a given set
of walks is the product of $w$ weights over all edges occupied by the
walks multiplied by the product of the $v$ weights for each of the
initial sites occupied. The generating function $\ZNWW_{t}(\bxi\pto
\bxf)$, of
$N$ walks of length $t$ starting at $\by=\bxi$ in column $m=0$ and
finishing at $\by=\bxf$ in column $m=t$ is the sum of these weights
over all sets of walks connecting $\bxi$ and $\bxf$:
\begin{equation}\label{part-def}
\ZNWW_{t}(\bxi\pto \bxf) = \sum_{\mathcal{Y}}\;
\prod_{j=1}^N v(y_j(0))\prod_{m=1}^t  w(y_j(m-1),y_j(m)) 
\end{equation}
where $y_j(m)$ is the position of the $j^{th}$ walk in column $m$ and
the set $\mathcal{Y}$ is given by
\begin{eqnarray}
\mathcal{Y} & = &\{ y_j(m)|1\leq j \leq N, 0 \leq m \leq t,1\leq y_1(m)< y_2(m) < \cdots < y_N(m) \leq L,
\nonumber \\
	& &  y_j(m)=y_j(m-1)\pm 1 \mbox{ and
} y_j(0) = y_j^i,\; y_j(t) = y_j^f.
\}
\end{eqnarray}
With homogeneous weights away from the boundaries
but extra weights at the boundaries the associated six-vertex model
has been considered in \cite{owczarek89}. However, we note
that in \cite{owczarek89}, and in most other studies of the
six-vertex model, only such properties of the model are calculated
that are averages over all numbers of walks, $N$. Here in contrast we
are considering the generating function for a fixed number of walks, $N$, of a
fixed finite length $t$. We will use the transfer matrix method from
statistical mechanics to find this generating function in terms of a
determinant of one-walk generating functions.

\begin{figure}[ht]
\begin{center}
 \figswitch{file=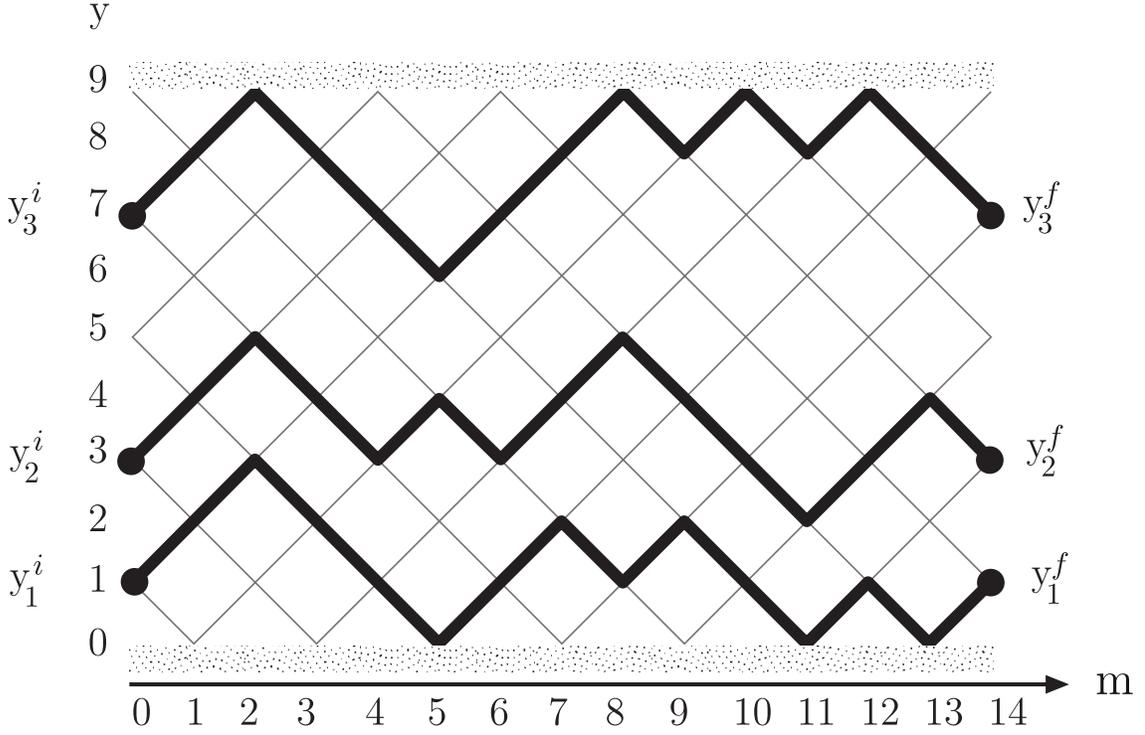}{width=15cm}
%
%
 \end{center}
\caption{\it   Three non-intersecting directed walks of length $t=14$
in a strip of width $L=9$. The variables $m$, $y_j^i$ and $y_j^f$ shown. The
walk closest to the lower wall has weight
$v(1)w(1,2)w(2,3)w(3,2)\cdots w(1,0)w(0,1)$ with $y^i_1=1$ and
$y^f_1=1$.  }
\label{fig:walks} 
\end{figure}

Precisely, we will show that the partition or generating function,
$\ZNWW_{t}(\byi\pto\byf)$, for non-intersecting walk configurations of
$N$ walks in which the $j^{th}$ walk starts at $y_j^i$ and arrives at
$y_j^f$ after $t$ steps is given by the following determinant:
\begin{equation} 
\ZNWW_{t}({\bf y}^i\pto {\bf y}^f)=
\left|
\begin{array}{cccc}
	\ZSWW_{t}(y_1^i\pto y_1^f)&\ZSWW_{t}(y_1^i\pto y_2^f)&
	   \dots&\ZSWW_{t}(y_1^i\pto y_N^f)\\
	\ZSWW_{t}(y_2^i\pto y_1^f)&\ZSWW_{t}(y_2^i\pto y_2^f)&
	   \dots&\ZSWW_{t}(y_2^i\pto y_N^f)\\
	.&.&.&.\\
	.&.&.&.\\
	.&.&.&.\\
	\ZSWW_{t}(y_N^i\pto y_1^f)&\ZSWW_{t}(y_N^i\pto y_2^f)&
	\dots&\ZSWW_{t}(y_N^i\pto y_N^f)
\end{array}
\right|
\label{eq:det2}
\end{equation}
where $\ZSWW_{t}(y_j^i\pto y_k^f)$ is the generating
function for configurations of a single walk starting at $y_j^i$ and
ending at $y_k^f$ in a strip of width $L$. This is then a
generalisation of the `master formulae' of Fisher (equation (5.9) of
\cite{fisher84}) and of Forrester (equation (4) of \cite{forrester89})
where unweighted non-intersecting walks are considered. Importantly,
this determinantal result is precisely that given by the very general
Gessel-Viennot Theorem \cite{gessel85,gessel89} for this problem: the
walks considered on this lattice with the generalised weights
considered satisfy the conditions of the theorem (see
\cite{guttmann98a} for a discussion of a one wall case with
homogeneous weights). 

\section{From Bethe Ansatz to determinant}
\label{sec:N-walkers}
\subsection{Transfer matrix formulation}

The generating function of our walk problem in a strip can be formulated as the
matrix element of a product of so-called transfer matrices (see
below for their precise definition).  This `transfer matrix' contains
the weights of all the possible edge configurations of $N$ walks of
two adjacent columns of sites. As such it is related to an invariant
subspace of the full six-vertex transfer matrix
\cite{guttmann98a}. In other words the non-zero elements of the
six-vertex transfer matrix occur in diagonal blocks each of which is
an ``$N$-walk transfer matrix'' for some $N$. There are of course
several ways of setting up a strip transfer matrix.  The most common
are the site-to-site and the edge-to-edge matrices which add either
edges or sites respectively to the walks.  In the case of the
six-vertex model the edge-to-edge matrix is usually chosen.  Our
choice of vertex-to-vertex (we shall use vertex and site
interchangeably) matrix is determined by the requirements that the
corresponding eigenvalue problem be as simple as possible while
ensuring the non-intersecting constraint can be easily implemented.

The calculation of the generating function using the transfer matrix
formulation then hinges on the spectral decomposition of the matrix.
The matrices we diagonalise will be matrices that add two-steps to the
evaluation of the generating function at a time. However the
generating function is initially constructed in terms of one-step
transfer matrices as these are useful since the non-intersecting condition
is explicit.

\begin{defin}\label{tm-def}
Let $y\in\cS{e}{L} $ and $y'\in\cS{o}{L}$. For $N=1$ the one-step 
transfer matrices are defined as
\begin{subequations}\label{eq:t-one}
	\begin{align}
	\left(\Teo\right)_{y,y'}&=
		\begin{cases}
			0         &\text{if $|y-y'|>1$}\\
			 w(y,y')  &\text{if $|y-y'|=1$}
\end{cases}\label{eq:teo}\\
\intertext{and}
	\left(\Toe\right)_{y',y}&=
		\begin{cases}
			0        &\text{if $|y'-y|>1$}\\
			 w(y',y) &\text{if $|y-y'|=1$}.
\end{cases}\label{eq:toe}
\end{align} 
\end{subequations}
The $N$-walk transfer matrices for $N>1$ are constructed from
sub-matrices of a direct product of the above $N=1$  
matrices:
\begin{subequations}\label{eq:t-N}
	\begin{align}
	(\TNoe)_{\by,\by'}&
	=\left(\bigotimes_{i=1}^{N}\Toe\right)_{\by,\by'}\qquad \qquad
	\text{$\by\in\cU{o}{L}$ and $\by' \in \cU{e}{L}$}
	\label{eq:dptmtoe}\\
\intertext{and}
(\TNeo)_{\by',\by}&
	=\left(\bigotimes_{i=1}^{N}\Teo\right)_{\by',\by}\qquad \qquad
	\text{$\by' \in \cU{e}{L}$ and $\by\in\cU{o}{L}$}
	\label{eq:dptmteo}
\end{align}
\end{subequations}
\end{defin}
Note that since we are only considering $L$ being odd we have that
$\TNoe$ and $\TNeo$ are square matrices and that in general
$\left(\Toe\right)_{y',y}\ne \left(\Teo\right)_{y,y'}$. More
importantly, the restriction of the row and column spaces of the
direct product eliminates the possibility of two walks arriving at the
same lattice point. Furthermore the condition that the one-walker
transfer matrix vanishes for $|y'-y|>1$ prevents the generation of
configurations in which pairs of walks ``cross'' \emph{without}
sharing a common lattice site (only nearest neighbour steps are
allowed in all cases). This ``non-crossing'' condition is
unnecessarily restrictive in the one-walk case. However, for $N>1$, if
further neighbour steps are allowed then, in general, it is \emph{not}
possible to use the Bethe Ansatz. This condition is the analogue of the
``non-crossing condition'' of the Gessel-Viennot Theorem.

The generating function $\ZNWW_{t}(\bxi\pto \bxf)$ of $N$
non-intersecting walks of length $t$ in a strip is related to
$\ZNWW_{t-1}(\bxi\pto \bf{y})$ by recurrence, the coefficients of which
are the elements of one of the two one-step transfer matrices defined
above.  This relationship is given by the following lemma.
\begin{lem}\label{part-tm}
The generating function $\ZNWW_{t}(\bxi\pto \bxf)$ is given for $t>0$, depending
on whether $\bxf \in \cU{e}{L}$ or not, by
\begin{equation}\label{transdef}
\ZNWW_{t}(\bxi\pto \bxf)=
\left\{
\begin{array}{l}
\sum_{\bx \in \cU{e}{L}}\ZNWW_{t-1}(\bxi \pto {\bf y})(\TNeo)_{{\bf y},\bxf}
\qquad\hbox{for}\qquad \bxf \in \cU{o}{L}\\
\\
\sum_{\bx \in \cU{o}{L}}\ZNWW_{t-1}(\bxi\pto {\bf y})(\TNoe)_{{\bf y},\bxf}
\qquad\hbox{for}\qquad \bxf \in \cU{e}{L}\; .
\end{array}
\right. 
\end{equation}
\end{lem}
\begin{proof}
A simple proof of this Lemma can be constructed using induction on $t$.
\end{proof}

Together with the initial condition
\begin{equation}
\ZNWW_{0}(\bxi\pto \bxf)= \d_{\bxi,\bxf} V(\bxi),
\end{equation}
where
\begin{equation}\label{factor}
V(\bxi)=\prod_{\a =1}^N v(y^i_\a),
\end{equation}
equation \eqref{transdef} determines $\ZNWW_{t}(\bxi\pto \bxf)$.

A simple corollary of this Lemma (again shown by induction) is that
the partition function can be written in terms of ``two-step''
transfer matrices.
\begin{cor}\label{two-steo-cor}
The generating function $\ZNWW_{t}(\bxi\pto \bxf)$ is given, depending 
on whether $t$ is even or odd, as
\begin{subequations}\label{part-tstrans}
\begin{align}\label{part-tstransa}
\ZNWW_{2r}(\bxi\pto \bxf)&=
\left\{
\begin{array}{ll}
V(\bxi)\left((\TNe)^{r}\right)_{\bxi,\bxf}
& \hbox{for}\quad \bxi \in \cU{e}{L} \hbox{ and } \bxf \in \cU{e}{L}\\
\\
V(\bxi)\left((\TNo)^{r}\right)_{\bxi,\bxf}
& \hbox{for}\quad \bxi \in \cU{o}{L} \hbox{ and } \bxf \in \cU{o}{L}
\end{array}
\right.\\
\intertext{or } 
\label{part-tstransb}
\ZNWW_{2r+1}(\bxi\pto \bxf)&=
\left\{
\begin{array}{ll}
V(\bxi)\sum_{\bx \in
\cU{e}{L}}\left((\TNe)^{r}\right)_{\bxi,\bx}(\TNeo)_{{\bf y},\bxf} 
& \hbox{for}\quad \bxi \in \cU{e}{L} \hbox{ and } \bxf \in \cU{o}{L}\\
\\
V(\bxi)\sum_{\bx \in
\cU{o}{L}}\left((\TNo)^{r}\right)_{\bxi,\bx}(\TNoe)_{{\bf y},\bxf} 
& \hbox{for}\quad \bxi \in \cU{o}{L} \hbox{ and } \bxf \in \cU{e}{L}
\end{array}
\right.
\end{align}
\end{subequations}
for $r =0,1,2,\dots$, where the two-step transfer matrices are defined as
\begin{subequations}\label{eq:twostep}
\begin{align}
	\TNe &=\TNeo\, \TNoe \label{eq:twostepa}\\
	\TNo &=\TNoe \,\TNeo.  \label{eq:twostepb}
\end{align}
\end{subequations}
\end{cor}
The two-step transfer matrices correspond to adding \emph{two} steps to the
paths.  We show below that the generating function $\ZNWW_{t}(\bxi\pto
\bxf)$ may in fact be expressed in terms of the eigenvalues and
eigenvectors of the two-step matrices. In general the two-step
matrices are \emph{not} symmetric and hence one has to consider left
and right eigenvectors.  The $N$-walk eigenvectors will be constructed
from the one-walk eigenvectors via the Bethe Ansatz.  Used in this
context the Ansatz expresses the components of the $N$-walk
eigenvectors as a determinant of the components of one-walk
eigenvectors (see (\ref{eq:Bethe}) below).

\subsection{From transfer matrices to determinants}

This section contains our main results in the form of two theorems
that state under what conditions the $N$-walk generating function can
be written as the determinant (\ref{eq:det2}). In summary, our
theorems proven below show that the \emph{equivalence} of the Bethe
Ansatz in the form of equation \eqref{eq:Bethe} and the result of the
Gessel-Viennot Theorem in the form (\ref{eq:stripdet}) rests on
showing that for any given problem the Ansatz is sufficiently good to
provide a spanning (or complete) set of eigenvectors.  This in turn
depends only on the completeness of the one-walk eigenvectors
(conditions of Lemma \ref{lem:cyclic}).

In this paper it is not our purpose to prove that all choices of the
functions $w(y,y')$ which make up the elements of the one-walk
transfer matrices allow for the conditions of our theorems to be
satisfied. However, these theorems have been subsequently used to
analyse some interesting cases in detail \cite{brak98e} where as a
consequence we show that an arbitrary boundary weight on one boundary
(and special weight on the other) with homogeneous weights otherwise
satisfies the conditions of the theorems. Also, we point out that
since the Gessel-Viennot determinant holds for any function $w(y,y')$
it is almost certainly true that our conditions are satisfied always.
It is however our main purpose here to demonstrate how the Bethe
Ansatz gives the solution of the $N$-walk problem \emph{given} the
solution of the one-walk problem.

The first theorem states the conditions under which the $N$-walk
transfer matrices can be diagonalised using a Bethe Ansatz: 
the major condition is
that the one-walk transfer matrix problem can be solved --- see
Lemma \ref{lem:cyclic}. The second theorem basically states that if the Bethe
Ansatz gives a complete set of eigenvectors for the $N$-walk problem for any
$N$, then the $N$-walk generating function is a determinant of
one-walk generating functions.

To begin we set up the conditions that define a solution of the
one-walk problem: these will be the conditions our two theorems
require. 

\begin{lem}
\label{lem:cyclic}
Suppose that there exist linearly independent sets of column vectors
$\{\evec{o}{R}{k}\}_{k\in\cK_{1}}$ and 
  $\{\evec{e}{R}{k}\}_{k\in\cK_{1}}$, where $\cK_{1}$ is 
  some index set, which satisfy
\begin{equation}\label{eq:cyclic-rR}
		  \Teo\evec{o}{R}{k}=\l_{k}\evec{e}{R}{k}
		  \qquad\hbox{and}\qquad
		  \Toe\evec{e}{R}{k}=\l_{k}\evec{o}{R}{k}
\end{equation}
with $\l_{k} \in \mathbb{C}$, and which span the column spaces of
$\Teo$ and $\Toe$ respectively (in which case they are said to be
complete).  Further let $\Te$ and $\To$ be defined by
\eqref{eq:twostep} then
\begin{itemize}
\item[(i)]
$\evec{o}{R}{k}$ and $\evec{e}{R}{k}$ are right 
  eigenvectors of $\To$ and $\Te$ respectively with 
  eigenvalue $\l_{k}^{2}$.

\item [(ii)] 
corresponding sets $\{\evec{o}{L}{k}\}_{k\in\cK_{1}}$ and
  $\{\evec{e}{L}{k}\}_{k\in\cK_{1}}$ of row vectors may be found such
  that
\begin{equation}\label{oc}
\evecb{p}{L}{k}\cdot\evec{p}{R}{k'}=\d_{k,k'} \qquad \hbox{and}
\qquad \sum_{k\in\cK_{1}}\evec{p}{R}{k}(y)\evecb{p}{L}{k}(y')=\d_{y,y'}
\end{equation}
for each $p\in\{e,o\}$, where the $*$ denotes complex
conjugation. Note that the vectors of \eqref{oc} have components
indexed by $y$.
\item[(iii)] the row vectors of (ii) satisfy
\begin{equation}\label{eq:cyclic-rL}
\evec{o}{L}{k}\Toe=\l_{k}\evec{e}{L}{k}
		  \qquad\hbox{and}\qquad
\evec{e}{L}{k}\Teo=\l_{k}\evec{o}{L}{k}
\end{equation}
  and also $\evec{o}{L}{k}$ and $\evec{e}{L}{k}$ are left 
  eigenvectors of $\To$ and $\Te$ respectively with eigenvalue 
  $\l_{k}^{2}$.
\end{itemize}
\end{lem}
Note that since the vectors in \eqref{eq:cyclic-rR} span the space the
cardinality of the index set $\cK_{1}$ is $(L+1)/2$.  The proof of the
lemma is elementary linear algebra and we omit it. Notice that if
$\l_k$ is a solution of \eqref{eq:cyclic-rR} then so is $-\l_k$ with
vector $\evec{e}{R}{k}$ replaced by $-\evec{e}{R}{k}$. These vectors
are clearly not independent and normally sufficient independent
vectors to form a spanning set are obtained by taking only the
positive values of $\l_k$.

From the above left and right one-walk vectors 
we now construct the $N$-walk vectors and hence eigenvectors of $\TNe$ and $\TNo$.
\begin{thm}\label{th:det} 
	 Let $\TNe$ and $\TNo$, $N>1$, be given by equations 
	 \eqref{eq:twostep}. By imposing an arbitrary ordering on the
elements of $\cK_{1}$ define  
         \begin{equation}
          \cK_{N}=\{\bk=(k_1,k_2,\dots k_N)|k_i\in\cK_{1}\mbox{ and } 
k_1<k_2<\dots < k_N\}
         \end{equation}
         and 
         \begin{equation}
		\Lam_{\bk}=\prod_{\a=1}^{N}\l_{k_{\a}}.
		\label{eq:Neval}
	 \end{equation} 
         (a) If for  $C\in\{L,R\}$ and  
	 $p\in\{e,o\}$, $\{\evec{p}{C}{k}\}_{k\in\cK_{1}}$ satisfy the
conditions of Lemma \ref{lem:cyclic} then the vectors
$\{\Evec{p}{C}{\bk}\}_{\bk\in\cK_{N}}$ given by the Bethe Ansatz,
	\begin{align}
		\Evec{p}{C}{\bk}(\by)=&\sum_{\s\in P_{N}}\e_{\s} 
		\prod_{\a=1}^{N}\evec{p}{C}{k_{\s_{\a}}}(y_{\a})
                =\sum_{\s\in P_{N}}\e_{\s} 
		\prod_{\a=1}^{N}\evec{p}{C}{k_{{\a}}}(y_{\s_\a})
		\qquad\text{$\bx\in\cU{p}{L}$},
		\label{eq:Bethe} 
	\end{align}
	where $P_{N}$ is the set of $N!$ permutations of
$\{1,2,\ldots,N\}$, $\s = (\s_{1}, \s_{2}, \ldots, \s_{N}) \in P_N$
and $\e_{\s}$ is the signature of the permutation $\s$,
        satisfy 
        \begin{equation}\label{eq:cyclic-RR}
		  \TNoe\Evec{e}{R}{\bk}=\Lam_{\bk}\Evec{o}{R}{\bk}
		  \qquad\hbox{and}\qquad
		  \TNeo\Evec{o}{R}{\bk}=\Lam_{\bk}\Evec{e}{R}{\bk}.
	\end{equation}
        (b) Moreover the conclusions of parts (i), (ii) and (iii) of Lemma
        \ref{lem:cyclic} hold with $\evec{p}{C}{k}$ replaced by
        $\Evec{p}{C}{\bk}$, $\Teo$ and $\Toe$ replaced by
        $\TNeo$ and $\TNoe$, $\cK_{1}$ replaced by $\cK_{N}$, and
        $\l_k$ replaced by $\Lam_{\bk}$.

\end{thm}
\noindent The proofs of part of this theorem and Theorem \ref{thm:specdecom} require the following result.
\begin{prop}
For $\bk\in\,\cK_{N}$ and $\by \in\,\cU{p}{L}$ let
\begin{equation}
\Phi_\bk(\bx) = \sum_{\s\in P_{N}}\e_{\s} 
		\prod_{\a=1}^{N}\phi_{k_{{\s_\a}}}(y_{\a})\qquad\hbox{and}\qquad
\Psi_\bk(\bx) = \sum_{\s\in P_{N}}\e_{\s} 
		\prod_{\a=1}^{N}\psi_{k_{\s_\a}}(y_{\a}).
\end{equation}
Also let
\begin{equation}
f(\bk)=\prod_{\a=1}^{N}f(k_{\a})
\end{equation}
then
\begin{equation}\label{ksum}
\sum_{\bk\,\in\,\cK_{N}}f(\bk)\Phi_\bk(\bx)\Psi_\bk(\bx')=\sum_{\s\in P_{N}}\e_{\s}
\prod_{\a=1}^{N}\left(\sum_{k_\a\in\,\cK_1 }f(k_\a)\phi_{k_{{\a}}}
(y_{\a}) \psi_{k_{\a}}(y_{\s_\a}')\right)
\end{equation}
and
\begin{equation}\label{ysum}
\sum_{\by\,\in\,\cU{p}{L}}\Phi_\bk(\bx)\Psi_{\bk'}(\bx)=\sum_{\s\in P_{N}}\e_{\s}
\prod_{\a=1}^{N}\left(\sum_{y_\a\,\in\cS{p}{L}}\phi_{k_{{\a}}}
(y_{\a}) \psi_{k_{\s_\a}'}(y_{\a})\right)
\end{equation}
\end{prop}
\begin{proof}
\begin{align}
\sum_{\bk\,\in\,\cK_{N}}f(\bk)\Phi_\bk(\bx)\Psi_\bk(\bx')=&\sum_{\s\in P_{N}}\e_{\s}
\sum_{\bk\,\in\,\cK_{N}}\sum_{\tau\in P_{N}}\e_{\tau}\prod_{\a=1}^{N}
f(k_\a)\phi_{k_{{\tau_\a}}}(y_{\a})\psi_{k_{{\s_\a}}}(y_{\a}')\\
=&\sum_{\s'\in P_{N}}\e_{\s'}
\sum_{\bk\,\in\,\cK_{N}}\sum_{\tau\in P_{N}}\prod_{\a=1}^{N}
f(k_\a)\phi_{k_{{\tau_\a}}}(y_{\a})\psi_{k_{{\tau_\a}}}(y_{\s'_\a}')
\end{align}
The double sum over permutations, $\tau$ and $\bk\in\cK_N$ is equivalent to summing each
$k_\a$ independently over $\cK_1$ (terms for which two or more components of $\bk$ are
equal make zero contribution) and the first result follows. The second result follows
in the same way by interchanging the roles of $k$ and $y$.
\end{proof}

\begin{proof}({\it of Theorem \ref{th:det}}$\,$) 
(a) We first obtain 
the cyclic property \eqref{eq:cyclic-RR} as follows.
\begin{subequations}\label{eq:crux}
\begin{align}
\left(\TNoe \Evec{e}{R}{\bk}\right)_\bx
&= \sum_{\by'\in\cU{e}{L}}\sum_{\s\in P_{N}}\e_{\s}
	\left(\TNoe\right)_{\by,\by'}
	\prod_{\a=1}^{N}\evec{e}{R}{k_{\s_{\a}}}(y'_{\a}) 
	\label{eq:crux1}\\
&= \sum_{\s\in P_{N}}\e_{\s}
	\sum_{\by'\in\cU{e}{L}}
	\prod_{\a=1}^{N}\left(\Toe\right)_{y_{\a},y'_{\a}}
	\evec{e}{R}{k_{\s_{\a}}}(y'_{\a})
	\qquad \left(\hbox{using}\quad\eqref{eq:t-N}\right)\label{eq:crux2} \\
&=\sum_{\s\in P_{N}}\e_{\s}\left[\sum_{y'_{1}\in\cS{e}{L}}
	\left(\Toe\right)_{y_{1},y'_{1}}\evec{e}{R}{k_{\s_{1}}}(y'_{1})
	\right] \ldots 
	\left[\sum_{y'_{N}\in\cS{e}{L}}
	\left(\Toe\right)_{y_{N},y'_{N}}
	\evec{e}{R}{k_{\s_{N}}}(y'_{N})
	\right]\label{eq:crux3}\\
&=\sum_{\s\in P_{N}}\e_{\s} \,
	\left[\l_{k_{\s_{1}}} \evec{o}{R}{k_{\s_{1}}}(y_{1})
	\right]
	\ldots  
	\left[\l_{k_{\s_{N}}}\evec{o}{R}{k_{\s_{N}}}(y_{N})
	\right] \qquad \left(\hbox{using}\quad\eqref{eq:cyclic-rR}\right)
	\label{eq:crux4}\\
&=  \Lam_{\bk}\, \Evec{o}{R}{\bk}(\by)
 	\end{align}
\end{subequations}
The critical step, and the whole reason for 
introducing the Bethe Ansatz, is to enable one to go from 
the restricted sums of \eqref{eq:crux2} to the unrestricted 
sums in \eqref{eq:crux3}.  This is justified for two 
reasons, 
\begin{enumerate}
\item since $\Evec{e}{R}{\bk}(\by')$ is a determinant, 
if any of the $y_{\a}$'s are equal then $\Evec{e}{R}{\bk}=0$ 
-- this allows the restriction $y'_{1}<y'_{2}\ldots<y'_{N}$ 
on the sum to be relaxed to $y'_{1}\le y'_{2}\ldots\le y'_{N}$
\item  the $y_\a$ are in strictly increasing order combined with 
the fact that the matrix elements of $\Toe$,
are only non-zero if $|y_{\a}-y'_{\a}|\le 1$ allows the 
restriction on the sum  to be removed altogether.  
\end{enumerate}
The second part of \eqref{eq:cyclic-RR} 
follows mutatis mutandis.

(b) (i) The vector $\Evec{e}{R}{\bk}$ is a right eigenvector of $\TNe$
with eigenvalue $\Lam_{\bk}^{2}$, since
\[
\TNe \Evec{e}{R}{\bk}=\TNeo\TNoe \Evec{e}{R}{\bk}=
\TNeo \Lam_{\bk} \Evec{o}{R}{\bk}= \Lam_{\bk}^{2} \Evec{e}{R}{\bk}
\]
which follows from \eqref{eq:cyclic-RR}. Similarly $\Evec{o}{R}{\bk} $ 
is a right eigenvector of $\TNo$ with eigenvalue $\Lam_{\bk}^{2}$.


\noindent (ii)  
Let us start by deriving the first result of part (ii), namely the
orthogonality and normalisation condition. Using \eqref{ysum}, for
$\bk,\bk'\in\cK_N$
\begin{align}
\sum_{\by \in \cU{p}{L} }\Evecb{p}{L}{\bk} (\bx)
	\Evec{p}{R}{\bk'} (\bx)\notag
=& \sum_{\s\in P_{N}}\e_{\s}
\prod_{\a=1}^{N}\left( \sum_{y_\a\,\in\cS{p}{L}}\evecb{p}{L}{k_\a}
(y_{\a}) \evec{p}{R}{k_{\s_\a}'}(y_{\a}) \right)\notag\\
=& \sum_{\s\in P_{N}}\e_{\s}\prod_{\a=1}^{N}\d_{k_\a,k_{\s_\a}'}
\qquad\left(\hbox{using}\quad\eqref{oc}\right)\notag\\
=&\prod_{\a=1}^{N}\d_{k_\a,k_{\a}'}\notag
\end{align}
since the components of $\bk$ and $\bk'$ are in the same order only the identity
permutation gives a non-zero delta function product. Thus
\begin{equation}\label{ortho}
\Evecb{p}{L}{\bk} \cdot \Evec{p}{R}{\bk} = \d_{\bk,\bk'}.
\end{equation}

Our derivation of the second part of (ii), namely the ``completeness
condition'', closely parallels the derivation of the orthogonality and
normalisation condition. Using
\eqref{ksum}, for $\bx,\bx'\in\cU{p}{L}$
\begin{align}
\sum_{\bk \in \cK_N }\Evec{p}{R}{\bk}(\bx)\Evecb{p}{L}{\bk}(\bx')\notag
=& \sum_{\s\in P_{N}}\e_{\s}
\prod_{\a=1}^{N}\left( \sum_{k_\a\,\in\cK_1}\evec{p}{R}{k_{\a}}(y_{\a})\evecb{p}{L}{k_\a}
(y_{\s_\a}')  \right)\notag\\
=& \sum_{\s\in P_{N}}\e_{\s}\prod_{\a=1}^{N}\d_{y_\a,y_{\s_\a}'}
\qquad\left(\hbox{using}\quad\eqref{oc}\right)\notag\\
=&\prod_{\a=1}^{N}\d_{y_\a,y_{\a}'}\notag
\end{align}
so
\begin{equation}\label{comp}
\sum_{\bk \in \cK_N }\Evec{p}{R}{\bk}(\bx)\Evecb{p}{L}{\bk}(\bx')= \d_{\bx,\bx'}
\end{equation}
Notice that $|\cK_N|= {{\half(L+1)}\choose N}$ which is the row (and column) space
dimension, as it should be for completeness.

\noindent (iii) Using basic linear algebra gives
\begin{equation}\label{eq:cyclic-RR-new}
		  \Evec{o}{L}{\bk}\TNoe=\Evec{e}{L}{\bk}\Lam_{\bk}
		  \qquad\qquad\qquad
		  \Evec{e}{L}{\bk}\TNeo=\Evec{o}{L}{\bk}\Lam_{\bk}
	\end{equation}
and
\begin{equation}\label{eq:cyclic-RR-new-t}
\Evec{e}{L}{\bk}\TNe = \Evec{e}{L}{\bk}\Lam_{\bk}^{2}\qquad\qquad\qquad
\Evec{o}{L}{\bk}\TNo = \Evec{o}{L}{\bk}\Lam_{\bk}^{2}
\end{equation}
\end{proof}
 
\begin{lem}\label{lem:part}
If the conditions of Theorem \ref{th:det} hold then
\begin{equation}\label{eq:sdN}
		\ZNWW_{t}(\bxi\pto \bxf)=
		V(\bx^i)
		\sum_{\bk\in\cbK}\,\,
		   \Evec{p'}{R}{\bk} (\bxi) \Lam^{t}_{\bk}
			\Evecb{p}{L}{\bk} (\bxf) 
\qquad\bxi\in\cU{p'}{L}\text{and}\quad\bxf\in\cU{p}{L}
\end{equation}
where if $t$ is even, $p'=p$ but otherwise $p$ and $p'$ are of opposite parity.
\end{lem}
\begin{proof}
If the conditions of Theorem \ref{th:det} hold then we have that
\eqref{ortho}, \eqref{comp},
\eqref{eq:cyclic-RR-new} and  \eqref{eq:cyclic-RR-new-t} are valid.  
Using \eqref{eq:cyclic-RR-new} and \eqref{comp} it follows that
\begin{equation}
(\,\TNoe)_{\bx,\bx'}=\sum_{\bk\in\cbK}\Lam_{\bk}\Evec{o}{R}{\bk} (\bx)\Evecb{e}{L}{\bk} (\bx')
\qquad\hbox{and}\qquad(\,\TNeo)_{\bx,\bx'}=
\sum_{\bk\in\cbK}\Lam_{\bk}\Evec{e}{R}{\bk} (\bx)\Evecb{o}{L}{\bk}
(\bx').
\end{equation}
Also, using \eqref{eq:cyclic-RR-new-t} and \eqref{comp} it follows that
\begin{equation}(\,\TNe)_{\bx,\bx'}=\sum_{\bk\in\cbK}\Lam_{\bk}^2\Evec{e}{R}{\bk} (\bx)\Evecb{e}{L}{\bk} (\bx')
\qquad\hbox{and}\qquad(\,\TNo)_{\bx,\bx'}=
\sum_{\bk\in\cbK}\Lam_\bk^2\Evec{o}{R}{\bk} (\bx)\Evecb{o}{L}{\bk}
(\bx').
\end{equation}
Substituting these into \eqref{part-tstrans} and using \eqref{ortho}
gives the result immediately.
%
\end{proof}

\begin{thm} 
	\label{thm:specdecom}
	 If the conditions of Theorem \ref{th:det} hold then
	 \begin{equation}
		\ZNWW_{t}(\bx^i\pto {\bf y}^f) =
		\det|| \ZSWW_{t}(y^i_\a\pto {y}^f_{\beta})
		||_{\a,\beta=1\ldots N}.
		\label{eq:stripdet}
	\end{equation}
\end{thm}
\begin{proof}
Using \eqref{factor}, \eqref{ksum} and \eqref{eq:sdN} (which follows
from \eqref{oc} by Lemma \ref{lem:part})  
\begin{align}
\ZNWW_{t}(\by^i\pto \by^f) 
&=\sum_{\s\in P_{N}}\e_{\s}
\prod_{\a=1}^{N}\left(v(y_\a^i)\sum_{k_{\a}\in\,\cK_1}\lambda_{k_{\a}}^t\evec{p'}{R}{k_{\a}}
(y^i_{\a}) \evecb{p}{L}{k_{\a}} (y^f_{\s_\a})\right)
		\notag\\
&= \sum_{\s\in P_N}
		\e_\s  \prod_{\a=1}^N 
		\ZSWW_{t}(y^i_\a
		\pto {y}^f_{\s_{\a}}).
		\label{eq:npart6}
\end{align}
which is an expansion of the required determinant. 
\end{proof}

Finally, we point out that the result of Theorem \ref{thm:specdecom}
is also the conclusion of the Gessel-Viennot Theorem.

\subsection{One wall and no wall geometries}

When $L > (t-|y^i_N-y^f_N|)/2+\max(y^i_N,y^f_N)$ the walk closest to
the wall at $y=L$ cannot touch it since the $N^{th}$ walk (of $t$
steps) needs at least $|y^i_N-y^f_N|$ steps to go from $y^i_N$ to
$y^f_N$ and any excursion close to the wall from the end point closest
to the wall ($\max(y^i_N,y^f_N)$) needs just as many steps to return
(hence the factor of $1/2$).  Hence, when this condition holds, the
strip generating function, $\ZSWW_{t}(y^i_\a \pto {y}^f_{\beta})$
becomes equal to the generating function for walks that are affected
by only one wall, $\ZSW_{t}(y^i_\a\pto {y}^f_{\beta})$.  Hence taking
the limit $L\to\infty$ gives the following corollary:
\begin{cor} For $\by^{i}\in\cU{p}{L}$ and $\by^{f}\in\cU{p'}{L}$, the
$N$-walk generating function with only one  wall at height $y=0$ is given by,
	\begin{equation}
		\ZNW_{t}(\by^i\pto \by^f)	= \det ||
		\ZSW_{t}(y^i_\a\pto {y}^f_{\beta})
		||_{\a,\beta=1\ldots N}\,. 
		\label{eq:semidet}
	\end{equation}	
\end{cor}

If we also condition the walk closest to the
wall at $y=0$ so that it cannot touch that wall we will end up with
the ``no boundary'' results.

\section*{Acknowledgements} 

Financial support from the Australian Research Council is gratefully
acknowledged by RB and ALO. JWE is grateful for financial support from
the Australian Research Council and for the kind hospitality provided
by the University of Melbourne during which time this research was
begun.



\end{document}